\documentclass[ 11pt,twoside]{article}
\usepackage{amsmath,amssymb,color}

   \usepackage{amsmath}
    \allowdisplaybreaks  

\def\BMO{{\mathrm{BMO}}}
\def\loc{{\mathrm{loc}}}

\def\Lip{{\mathrm{Lip}}}

\def\Osc{{\mathrm{Osc}}}
\def\dint{\displaystyle\int}

\usepackage{graphicx,amssymb,mathrsfs,amsmath,latexsym,amsfonts,titlesec,amscd,epsfig,cite,amsthm}
\usepackage{indentfirst} 

\newcommand{\upcite}[1]{\textsuperscript{\textsuperscript{\cite{#1}}}}
\numberwithin{equation}{section}
\theoremstyle{definition} 

\newtheorem{theorem}{\indent  Theorem}[section]
    \newtheorem{lemma}{\indent  Lemma} [section]
\newtheorem{corollary}{\indent Corollary}[section]

        \newtheorem{remark}{\indent  Remark}  

\newtheorem*{acknowledgments}{\indent Acknowledgments}
\makeatletter
\def\address#1#2{\begingroup
\noindent\parbox[t]{7.8cm}{%
\small{\scshape\ignorespaces#1}\par\vskip1ex
\noindent\small{\itshape E-mail address}%
\/: #2\par\vskip4ex}\hfill%
\endgroup}%
\makeatother
\title{\bf\Large\uppercase{
Weighted Estimates for Multilinear Commutators of Marcinkiewicz Integrals with Bounded Kernel}} 
\author{ \textsc{\normalsize Jianglong Wu$^{1}$ \ \ Qingguo Liu$^{2}$  } 
 }
\date{} 
\begin{document}

\maketitle
\noindent\footnote{

\rightline{ \scriptsize \sc Ukrainian Mathematical Journal,  Vol. {\bf\rm 65} No. X, 2013.}
}

\begin{minipage}[t]{15cm}  \setlength{\baselineskip}{1.0em}
\noindent  
 {\bf Abstract:} Let $\mu_{\Omega,\vec{b}}$ be the multilinear commutator generalized by $\mu_{\Omega}$, the $n$-dimensional Marcinkiewicz integral with the bounded kernel, and $b_{j}\in \Osc_{\exp L^{r_{j}}}(1\le j\le m)$. In this paper, the following weighted inequalities are proved for $\omega\in
A_{\infty}$ and $0<p<\infty$,
$$\|\mu_{\Omega}(f)\|_{L^{p}(\omega)}\leq C\|M(f)\|_{L^{p}(\omega)},
\ \ \|\mu_{\Omega,\vec{b}}(f)\|_{L^{p}(\omega)}\leq C\|M_{L(\log
L)^{1/r}}(f)\|_{L^{p}(\omega)}.$$ The weighted weak $L(\log
L)^{1/r}$ -type estimate is also established when $p=1$ and
$\omega\in A_{1}$.
\\

 {\bf Keywords:}  \  Marcinkiewicz integral; multilinear commutator;  sharp function; $A_{p}$
 weight
 \\

 { \bf MR(2000) Subject Classification:}  \ 42B20; 42B25; 42B30  \ \

\end{minipage}

\section{Introduction and Main Results}

Suppose that $S^{n-1}$ is the unit sphere in $\mathbf{R}^{n}(n\geq
2)$ equipped with the normalized Lebesgue measure
$\mathrm{d}\sigma$. Let $\Omega\in L^{1}(S^{n-1})$ be a homogeneous  function of degree zero which satisfies  the cancellation  condition
\begin{equation} \label{equ:3.1}
\int_{S^{n-1}}\Omega(x')\mathrm{d}x'=0,
\end{equation}
where $x'=x/|x|~(\forall x\ne 0)$.

The n-dimensional Marcinkiewicz integral corresponding to the
Littlewood-Paley $g$-function introduced by Stein\upcite{S2} is
defined by
$\mu_{\Omega}(f)(x)=\big(\int_{0}^{\infty}|F_{\Omega,t}(f)(x)|^{2}\frac{\mathrm{d}t}{t^{3}}\big)^{1/2},
$
where
$F_{\Omega,t}(f)(x)=\int_{|x-y|\leq t}\frac{\Omega(x-y)}{|x-y|^{n-1}}f(y)\mathrm{d}y.
$

As usual, we denote by $A_{p}~(1\leq p<\infty)$  the Muckenhoupt's
weights class. We denote $[\omega]_{A_{p}}$ as $A_{p}$ constant (see \cite{S1} Chapter V or \cite{G} Chapter 9 for details). Operators that map $L^{p}$ to  $L^{q}$ are called  of strong type $(p, q)$ and operators that map $L^{p}$ to  $L^{q,\infty}$ are called  of weak type $(p, q)$ ~(see \cite{G} p. 32). Let
\begin{displaymath}
 \log^{+}t=\max (\log t,0) =\left\{\begin{array}{ll}
\log t, &\hbox{when}\ t>1, \\
0, &\hbox{when}\ 0\le t\le 1,
\end{array}\right.
\end{displaymath}
where $\log t=\ln t$, and we denote by $L(\log L)$
the set of all $f$ with $\int_{\mathbf{R}^{n}} |f(x)|\log^{+}|f(x)| \mathrm{d}x <\infty $ (see \cite{S1} p.128; \cite{G} \S7.5.a ). Here and in what follows, $\|b\|_{\ast}$ denotes the  $\BMO$  norm of $b$ (see  \cite{G} Chapter 7 for details).

In 1958, Stein\upcite{S2} proved  that $\mu_{\Omega}$ is  of  strong type $(p, p)$ for $1<p\leq 2$ and  of  weak  type  $(1, 1)$  when $\Omega\in \Lip_{\alpha}$ $(0<\alpha\leq 1)$, that is, there is a constant $C>0$
such that
\begin{equation} \label{equ:3.2}
|\Omega(x')-\Omega(y')|\leq C|x'-y'|^{\alpha},\ \ \forall~ x',y'\in S^{n-1}.
\end{equation}

In 1990, Torchinsky and Wang\upcite{TW} studied the weighted $L^{p}$
-boundedness of $\mu_{\Omega}$ when $\Omega$ satisfies
(\ref{equ:3.1}) and (\ref{equ:3.2}). They also considered the
weighted $L^{p}$ -norm inequality for the commutator of the
Marcinkiewcz integral, which is defined by
$$\mu_{\Omega, b}^{m}(f)(x)=\Big(\int_{0}^{\infty}\Big|\int_{|x-y|\leq t}
\frac{(b(x)-b(y))^{m}\Omega(x-y)}{|x-y|^{n-1}}f(y)\mathrm{d}y\Big|^{2}
\frac{\mathrm{d}t}{t^{3}}\Big)^{1/2},\ \ m\in\mathbf{N}.
$$

In 2004, Ding, Lu and Zhang\upcite{DLZ} studied the weighted weak
$L(\log L)$-type estimates for $\mu_{\Omega, b}^{m}$, precisely,  if
$\omega\in A_{1},~b\in \BMO$, $\Omega$ satisfies (\ref{equ:3.1}) and (\ref{equ:3.2}), then,
 for all $\lambda>0$,   there exists a constant $C>0$, such that
\begin{equation*} \label{equ:3.3}
\omega(\{x\in\mathbf{R}^{n}: |\mu_{\Omega, b}^{m}(f)(x)|>\lambda\}) \leq C \int_{\mathbf{R}^{n}}\frac{|f(x)|}{\lambda}\bigg(1+\log^{+}\frac{|f(x)|}{\lambda}\bigg)^{m}\omega(x)
\mathrm{d}x .
\end{equation*}

In 2008, Zhang\upcite{Z2} studied the weighted boundedness for the multilinear commutator of Marcinkiewicz integral $\mu_{\Omega, \vec{b}}$ when $\Omega\in \Lip_{\alpha}(0<\alpha\leq 1)$, $0<p<\infty$ and $\omega\in A_{\infty}$ (see  \cite{G} \S9.3), and  established a weighted weak $L(\log
L)^{1/r}$-type estimate when $p=1$ and $\omega\in A_{1}$, where
$$\mu_{\Omega, \vec{b}}(f)(x)=\bigg(\int_{0}^{\infty}\bigg|\int_{|x-y|{\leq}t}
\frac{\Omega(x-y)}{|x-y|^{n-1}}\Big(\prod_{j=1}^{m}
\big(b_{j}(x)-b_{j}(y)\big)\Big)f(y)\mathrm{d}y\bigg|^{2}
\frac{\mathrm{d}t}{t^{3}}\bigg)^{1/2},\ \ m\in\mathbf{N}.
$$
And in 2012, Zhang, Wu and Liu\upcite{ZWL} establish the weighted weak $L(\log L)^{m}$-type estimate for $\mu_{\Omega, \vec{b}}$ when  $\Omega$ satisfies a kind of Dini conditions.

In 2004, Lee and Rim\upcite{LR} proved the $L^{p}$ boundedness for
$\mu_{\Omega}$ when there exist constants $C>0$ and $\rho>1$ such that
\begin{equation} \label{equ:3.4}
|\Omega(x')-\Omega(y')|\leq\frac{C}{\bigg(\log\frac{1}{|x'-y'|}\bigg)^{\rho}}
\end{equation}
holds uniformly in $x', y'\in S^{n-1}$, and  $\Omega\in L^{\infty}(S^{n-1})$ be a homogeneous function of degree zero with cancellation property (\ref{equ:3.1}).
In 2005, Ding\upcite{D} studied the weak $(1,1)$ type estimate when
$\rho>2$ and $\Omega$ satisfies (\ref{equ:3.1}) and (\ref{equ:3.4}).

In the following, we will always assume that $\Omega\in L^{\infty}(S^{n-1})$ and satisfies (\ref{equ:3.1}) and (\ref{equ:3.4}), where $\rho>2$. Let $m$ be a positive integer. For
$\vec{b}=(b_{1}, b_{2},\cdots ,b_{m})$, $b_{j}\in \Osc_{\exp L^{r_{j}}}$,
$r_{j}\geq 1(1\leq j\leq m)$, we denote
\begin{equation} \label{equ:3.5}
\frac{1}{r}=\frac{1}{r_{1}}+\cdots+\frac{1}{r_{m}},\qquad \|\vec{b}\|=\prod_{j=1}^{m}\|b_{j}\|_{\Osc_{\exp L^{r_{j}}}}.
\end{equation}
For the definitions of $\Osc_{\exp L^{r}}$, $\|\cdot\|_{\Osc_{\exp
L^{r}}}$ and $M_{L(\log L)^{1/r}}$, see $\S$2.

Our results can be stated as follows.

\begin{theorem}\label{thm.1}\ \  Let $0<p<\infty$ and suppose that $\omega\in
A_{\infty}$. For $\rho>2$, $\Omega\in L^{\infty}(S^{n-1})$ is
homogeneous of degree zero and satisfies (\ref{equ:3.1}) and
(\ref{equ:3.4}). Then there is a positive constant $C$, such that
\begin{equation*} \label{equ:3.6}
\int_{\mathbf{R}^{n}}|\mu_{\Omega}(f)(x)|^{p}\omega(x)\mathrm{d}x
\leq C[\omega]^{p}_{A_{\infty}}\int_{\mathbf{R}^{n}}[M(f)(x)]^{p}\omega(x)\mathrm{d}x
\end{equation*}
for all bounded functions $f$ with compact support.
\end{theorem}

\begin{theorem}\label{thm.2}\ \ Let $0<p<\infty$, $\omega\in A_{\infty}$ and
$b_{j}\in \Osc_{\exp L^{r_{j}}}$, $r_{j}\geq 1~(1\leq j\leq m)$. $r$
and $\|\vec{b}\|$ be as in (\ref{equ:3.5}). For $\rho>2$, $\Omega\in
L^{\infty}(S^{n-1})$ is homogeneous of degree zero and satisfies
(\ref{equ:3.1}) and (\ref{equ:3.4}). Then there is a positive
constant $C$, such that
\begin{equation} \label{equ:3.7}
\int_{\mathbf{R}^{n}}|\mu_{\Omega,\vec{b}}(f)(x)|^{p}\omega(x)\mathrm{d}x
\leq C\|\vec{b}\|^{p}\int_{\mathrm{R}^{n}}[M_{L(\log L)^{1/r}}(f)(x)]^{p}\omega(x)\mathrm{d}x
\end{equation}
for all bounded functions $f$ with compact support.
\end{theorem}

Since $r_{j}\geq 1(j=1,2,\cdots,m)$, then $M_{L(\log L)^{1/r}}$  is
pointwise smaller than $M_{L(\log L)^{m}}$. Noting that $M_{L(\log
L)^{m}}$ is equivalent to $M^{m+1}$, the $m+1$ iterations of the
Hardy-Littlewood maximal operator $M$ (see (21) in \cite{P}), by using
the weighted $L^{p}$ -boundedness of $M$ again, from Theorem  \ref{thm.2}, we
have the following result.

\begin{corollary}\label{cor.1}\ \ Let $1<p<\infty$,  $\omega\in A_{p}$, $b_{j}\in \Osc_{\exp L ^{r_{j}}}$, $r_{j}\geq 1~(1\leq j\leq m)$, $r$ and $\|\vec{b}\|$ be as in (\ref{equ:3.5}). For $\rho>2$,
$\Omega\in L^{\infty}(S^{n-1})$ is homogeneous of degree zero and satisfies (\ref{equ:3.1}) and (\ref{equ:3.4}). Then there is a positive constant $C$, such that
\begin{equation*} \label{equ:3.8}
\int_{\mathbf{R}^{n}}|\mu_{\Omega,
\vec{b}}(f)(x)|^{p}\omega(x)\mathrm{d}x
\leq C\|\vec{b}\|^{p}\int_{\mathrm{R}^{n}}|f(x)|^{p}\omega(x)\mathrm{d}x
\end{equation*}
for all bounded functions $f$ with compact support.
\end{corollary}

\begin{theorem}\label{thm.3}\ \ Let $\omega\in A_{1}$,  $b_{j}\in \Osc_{\exp L
^{r_{j}}}$, $r_{j}\geq 1~(1\leq j\leq m)$, $r$ and $\|\vec{b}\|$ be
as in (\ref{equ:3.5}). For $\rho>2$, $\Omega\in L^{\infty}(S^{n-1})$
is homogeneous of degree zero and satisfies (\ref{equ:3.1}) and
(\ref{equ:3.4}). Let $\Phi(t)=t\log^{1/r}(\mathrm{e}+t)$. Then there
is a positive constant $C$, for all bounded functions
$f$ with compact support and all $\lambda>0$, such that
\begin{equation*} \label{equ:3.9}
\omega(\{x\in\mathbf{R}^{n}: \mu_{\Omega, \vec{b}}(f)(x)>\lambda\})
\leq C\int_{\mathbf{R}^{n}}\Phi\bigg(\frac{\|\vec{b}\||f(y)|}{\lambda}\bigg)\omega(y)\mathrm{d}y.
\end{equation*}
\end{theorem}

The remainder of the paper is organized as follows. In  \S2,  we will recall some notation  and known results we need, and establish the basic estimates for sharp functions.  In \S3 we prove Theorem \ref{thm.1} and \ref{thm.2}. In the last section, we prove Theorem \ref{thm.3}.

Throughout this paper, $C$ denotes a constant that is independent of
the main parameters involved but whose value may differ from line to
line. For any index $p \in [1, \infty]$, we denote by $p'$ its
conjugate index, namely, $1/p+1/p' = 1$. For $A \sim B$, we mean
that there is a constant $C > 0$ such that $C^{-1}B \le A \le CB$.

\section{Preliminaries and Estimates for Sharp Functions}

As usual, $M$ stands for the Hardy-Littlewood maximal operator. For
a ball $B$ in $\mathbf{R}^{n}$, denote  by
$f_{B}=|B|^{-1}\int_{B}f(y)\mathrm{d}y$. We need the following
variants of $M$ and the Fefferman-Stein's sharp function. For
$\delta>0$, define
$$M_{\delta}(f)(x)=\big[M(|f|^{\delta})(x)\big]^{1/\delta},\quad
M_{\delta}^{\sharp}(f)(x)=\big[M^{\sharp}(|f|^{\delta})(x)\big]^{1/\delta},
$$
where
$$M^{\sharp}(f)(x)=\sup_{B{\ni}x}\inf_{c}\frac1{|B|} \int_{B}|f(y)-c|\mathrm{d}y\approx \sup_{B{\ni} x} \frac1{|B|}\int_{B}|f(y)-f_{B}|\mathrm{d}y.
$$

The following relationships between $M^{\sharp}_{\delta}$ and $M_{\delta}$ which will be used is a version of the classical ones due to Fefferman and Stein (see \cite{S1} p.153).

\begin{lemma}\upcite{P,PT,P2} \label{lem.1}\ \ (a)  Let $\omega\in A_{\infty}$ and $\phi :  (0,\infty)\rightarrow(0,\infty)$ be doubling. Then there exists a positive constant $C$, depending upon the doubling condition of $\phi$, such that, for all $\lambda,~\delta>0$
$$\sup_{\lambda>0}\phi(\lambda)\omega(\{y\in \mathbf{R}^{n}: M_{\delta}(f)(y)>\lambda\})\leq C[\omega]_{A_{\infty}} \sup_{\lambda>0}\phi(\lambda)\omega(\{y\in\mathbf{R}^{n}:
M^{\sharp}_{\delta}(f)(y)>\lambda\}),$$
for every function $f$ such that the left-hand side is finite.

(b) Let $\omega\in A_{\infty}$ and $0<p,\delta<\infty$. Then there exists a positive constant $C$, depending upon $p$, such that $$\int_{\mathbf{R}^{n}}\big[M_{\delta}(f)(x)\big]^{p}\omega(x)\mathrm{d}x
{\leq}C[\omega]_{A_{\infty}}^{p}\int_{\mathbf{R^{n}}} \big[M_{\delta}^{\sharp}(f)(x)\big]^{p}\omega(x)\mathrm{d}x, $$
 for every function $f$ such that the left-hand side is finite.
\end{lemma}

A function $\Phi$ defined on $[0,\infty)$ is said to be a Young function, if $\Phi$ is a continuous, nonnegative, strictly increasing and convex function with $\Phi(0)=0$ and $\lim\limits_{t\rightarrow\infty}\Phi(t)=\infty$.
Define the $\Phi-$average of a function $f$ on a ball B by $$\|f\|_{\Phi,B}=\inf\bigg\{\lambda>0: \frac1{|B|}\int_{B} \Phi\bigg(\frac{|f(y)|}{\lambda}\bigg)\mathrm{d}y\leq 1\bigg\}. $$
The maximal operator $M_{\Phi}$ associated with the $\Phi-$average, $\|\cdot\|_{\Phi,B}$, is defined by
$$M_{\Phi}(f)(x)=\sup_{B\ni x}\|f\|_{\Phi,B},$$ where the supremum is taken over all the balls $B$ containing $x$.

When $\Phi(t)=t\log^{r}(\mathrm{e}+t)$, we denote
$\|\cdot\|_{\Phi,B}$ and $M_{\Phi}$ by $\|\cdot\|_{L(\log L)^{r},B}$
and $M_{L(\log L)^{r}}$, respectively. When
$\Phi(t)=\mathrm{e}^{t^{r}}-1$, we denote $\|\cdot\|_{\Phi,B}$ and
$M_{\Phi}$ by $\|\cdot\|_{\exp L^{r},B}$ and $M_{\exp L^{r}}$. If
$k\in\mathbf{N}$ then $M_{L(\log L)^{m}}\sim M^{m+1}$(see (21) of
\cite{P}).

We have the generalized H\"{o}lder's inequality as follows, for details and the more general cases see Lemma 2.3 in \cite{PT}.

\begin{lemma}\upcite{PT} \label{lem.2}\ \  Let $r_{1},\cdots,r_{m}\geq 1$ with
$1/r=1/r_{1}+\cdots+1/r_{m}$ and $B$ be a ball in $\mathbf{R}^{n}$.
Then there holds the generalized H\"older's inequality
$$\frac{1}{|B|}\int_{B}|f_{1}(x)\cdots f_{m}(x)g(x)|\mathrm{d}x\leq C\|f_{1}\|_{\exp L ^{r_{1}},B}\cdots\|f_{m}\|_{\exp L ^{r_{m}},B}\|g\|_{L(\log L)^{1/r},B}.$$
\end{lemma}

For $r\geq 1$, we say $f \in \Osc_{\exp L^{r}}$ if
$f{\in}L^{1}_{\mathrm{\loc}} (\mathbf{R}^{n})$ and $\|f\|_{\Osc_{\exp
L^{r}}}<\infty$, where
$$\|f\|_{\Osc_{\exp L^{r}}}=\sup_{B}\|f-f_{B}\|_{\exp L^{r},B},
$$
and the supremum is taken over all the balls $B\subset
\mathbf{R}^{n}$.

By John-Nirenberg theorem  (see \cite{S1} or \cite{JN}),  it is not
difficult to see that $\Osc_{{\exp}L^{1}} =\BMO(\mathbf{R}^{n})$ and
$\Osc_{\exp L^{r}}$ is contained properly in $\BMO(\mathbf{R}^{n})$
when $r>1$~(see \cite{HMY}). Furthermore, $\|b\|_{\ast}\leq
C\|b\|_{\Osc_{\exp L^{r}}}$ when $b \in \Osc_{\exp L^{r}}$ and $r\geq
1$~(see \cite{Z2}). For more information on Orlicz space see
\cite{RR}.

We will take the point of view of the vector-valued singular integral of Benedek, Calder\'{o}n and Panzone\upcite{BCP}. Let $\mathcal{H}$ be the Hilbert space defined by
$$\mathcal{H}=\bigg\{h:\; \|h\|_{\mathcal{H}}=\bigg(\int_{0}^{\infty} \frac{|h(t)|^{2}}{t^{3}}\mathrm{d}t\bigg)^{1/2}<\infty\bigg\}.
$$
For all  $x\in\mathbf{R}^n$ and $t>0$, let
$$F_{\Omega,\vec{b},t}(f)(x)=\int_{|x-y|{\leq}t}
\frac{\Omega(x-y)}{|x-y|^{n-1}}\Big(\prod_{j=1}^{m}
\big(b_{j}(x)-b_{j}(y)\big)\Big)f(y)\mathrm{d}y,\ \ m\in\mathbf{N}.
$$
Then for each fixed $x\in\mathbf{R}^n$, $F_{\Omega,t}(f)(x)$ and
$F_{\Omega,\vec{b},t}(f)(x)$ can be regarded as mapping from
 $[0,\infty)$ to $\mathcal{H}$, and
$$\mu_{\Omega}(f)(x)=\|F_{\Omega,t}(f)(x)\|_{\mathcal{H}},\ \
\mu_{\Omega,
\vec{b}}(f)(x)=\|F_{\Omega,\vec{b},t}(f)(x)\|_{\mathcal{H}}.
$$

The following pointwise estimates for the sharp function of $\mu$ come from \cite{L}.

\begin{lemma}\upcite{L}\label{lem.3}\ \   Let $0<\delta<1$, $f,~  \mu_{\Omega}(f)$ be both
locally integrable function. For $\rho>2$, $\Omega\in L^{\infty}(S^{n-1})$ is homogeneous of degree zero and satisfies (\ref{equ:3.1}) and (\ref{equ:3.4}). Then there is a positive constant $C$ ,
independent of $f$ and $x$, such that
$$M_{\delta}^{\sharp}(\mu_{\Omega}(f))(x)\leq CM(f)(x),\qquad a.e. \ x\in\mathbf{R}^{n}.$$
\end{lemma}

Some  ideas  for the proof of Lemma \ref{lem.3}  come from~\cite{DLZ}. For   details and the more information see Lemma 3.2.4 in \cite{L}.

For the multilinear commutators $\mu_{\Omega, \vec{b}}$,  there holds a similar piontwise estimate. To state it, we first introduce some  notations. For all $1{\leq} j{\leq} m$, we denote by $\mathscr{C}_{j}^{m}$  the family of all finite subsets $\sigma=\{\sigma(1),\cdots,\sigma(j)\}$ of $\{1,2,\cdots,m\}$ with $j$ different elements. For any $\sigma\in \mathscr{C}_{j}^{m}$ and $\vec{b}=(b_{1},\cdots,b_{m})$,
  we define $\sigma'=\{1,2,\cdots, m\}\setminus\sigma,$
  $\vec{b}_{\sigma}=(b_{\sigma(1)},\cdots,b_{\sigma(j)})$, and
$b_{\sigma}=b_{\sigma(1)}\cdots b_{\sigma(j)}$.  For any vector
$(r_{\sigma(1)},\cdots,r_{\sigma(j)})$ of $j$ positive numbers and
$1/r_{\sigma}=1/r_{\sigma(1)}+\cdots+1/r_{\sigma(j)}$, we write
\begin{equation} \label{equ:3.15}
\|\vec{b}_{\sigma}\|=\|\vec{b}_{\sigma}\|_{\Osc_{\exp L
^{r_{\sigma}}}}
=\|b_{\sigma(1)}\|_{\Osc_{\exp L ^{r_{\sigma(1)}}}}\cdots\|b_{\sigma(j)}\|_{\Osc_{\exp L ^{r_{\sigma(j)}}}}.
\end{equation}

For any
$\sigma=\{\sigma(1),\cdots,\sigma(j)\}\in\mathscr{C}_{j}^{m}$ and
$\vec{b}_{\sigma}=(b_{\sigma(1)},\cdots,b_{\sigma(j)})$, we write
$$F_{\Omega,\vec{b}_{\sigma},t}(f)(x)=\int_{|x-y|{\leq} t}
\frac{\Omega(x-y)}{|x-y|^{n-1}} \Big(\prod_{i=1}^{j}
\big(b_{\sigma(i)}(x)-b_{\sigma(i)}(y)\big)\Big)f(y) \mathrm{d}y
$$
and
$$\mu_{\Omega, \vec{b}_{\sigma}}(f)(x)
=\|F_{\Omega,\vec{b}_{\sigma},t}(f)(x) \|_{\mathcal{H}}.$$ If
$\sigma=\{1,\cdots,m\}$, then $\sigma'=\O$. We understand
$\mu_{\Omega, \vec{b}_{\sigma}}=\mu_{\Omega, \vec{b}}$ and
~$\mu_{\Omega, \vec{b}_{\sigma'}}=\mu_{\Omega}$.

\begin{lemma}\upcite{L}\label{lem.4}\ \ Let $r_{j}\geq 1$, $b_{j}\in \Osc_{\exp L
^{r_{j}}}(1\leq j\leq m)$, $r$ and $\|\vec{b}\|$ be as in
(\ref{equ:3.5}). For $\rho>2$, $\Omega\in L^{\infty}(S^{n-1})$ is
homogeneous of degree zero satisfying (\ref{equ:3.1}) and
(\ref{equ:3.4}), then for any $\delta$ and $\varepsilon$ with
$0<\delta<\varepsilon<1$, there is a constant $C>0$, depending only
on $\delta$ and $\varepsilon$, such that, for any bounded function
$f$ with compact support,
\begin{equation*} \label{equ:3.16}
M_{\delta}^{\sharp}(\mu_{\Omega, \vec{b}}(f))(x)\leq C\Big(
\|\vec{b}\|M_{L(\log L)^{1/r}}(f)(x) +\sum_{j=1}^{m} \sum_{\sigma\in
\mathscr{C}_{j}^{m}} \|\vec{b}_{\sigma}\|_{\Osc_{\exp L
^{r_{\sigma}}}}M_{\varepsilon}
\big(\mu_{\Omega, \vec{b}_{\sigma'}}(f)\big)(x)\Big).
\end{equation*}
\end{lemma}

Some  ideas  for the proof of Lemma \ref{lem.4}  come from~\cite{DLZ, P, PT,Z2}. For   details and the more information see Lemma 3.2.5 in \cite{L}.

\begin{remark}\ \ Noting that (\ref{equ:3.4}) is weaker than
$\Lip_{\alpha}(0<\alpha\leq 1)$ condition, the main results in this
paper improve the  main results in \cite{Z2}. And the Theorem \ref{thm.3}
is equivalent to the theorem 4.1.1 in \cite{SXF} when
$b_{1}=b_{2}=\cdots=b_{m}$.
\end{remark}

\section{Proof of Theorems 1.1 and 1.2}

The proof of Theorem \ref{thm.1} is similar as Theorem 1.1 in \cite{Z2}, So,
we omit the details  and only give the proof of Theorem \ref{thm.2} here. For brevity, we write
$$\|h(x)\|_{L^{p}(\omega)}= \Big( \int_{\mathbf{R}^{n}} |h(x)|^{p}\omega(x)\mathrm{d}x \Big)^{1/p}, \hbox{for} 0<p<\infty. $$

 \begin{proof} [\indent\bf Proof of Theorem \ref{thm.2}]\  Without loss of generality, we
assume
\begin{equation} \label{equ:3.34}
\int_{\mathbf{R}^{n}}[M_{L(\log L)^{1/r}}(f)(x)]^{p}\omega(x)\mathrm{d}x<\infty,
\end{equation}
since otherwise there is nothing to be proven. We divide the proof
into two cases.

Case I\quad  Suppose that $\omega$ and $b_{j}~(1\leq j\leq m)$ are all
bounded. Firstly, we take it for granted that, for all bounded
functions $f$ with compact supports,
\begin{equation} \label{equ:3.35}
\int_{\mathbf{R}^{n}}[M_{\delta}(\mu_{\Omega, \vec{b}}(f))(x)]^{p}\omega(x)\mathrm{d}x<\infty
\end{equation}
holds for $0<p<\infty$ and appropriate $\delta$ with $0<\delta<1$.

Under the assumption of (\ref{equ:3.35}), we will proceed the proof
by induction on $m$.  For $m=1$, $\vec{b}=b_{1}$, $\mu_{\Omega,
\vec{b}}=\mu_{\Omega, b_{1}}$. By Lemma 2.1 (b) and Lemma 2.4, for
$0<\delta<\varepsilon<1$, we have
\begin{equation}     \label{equ:3.10}
\begin{split}
\|\mu_{\Omega, b_{1}}(f)\|_{L^{p}(\omega)}&\leq \|M_{\delta} (\mu_{\Omega, b_{1}}(f))\|_{L^{p}(\omega)} \leq C \|M_{\delta}^{\sharp}(\mu_{\Omega, b_{1}}(f))\|_{L^{p}(\omega)}\\
&\leq C\|b_{1}\|_{\Osc_{\exp L ^{r_{1}}}}\Big(\|M_{L(\log L)^{1/r_{1}}}(f)\|_{L^{p}(\omega)}
+\|M_{\varepsilon}(\mu_{\Omega}(f))\|_{L^{p}(\omega)}\Big).
\end{split}
\end{equation}

 Since $\omega\in A_{\infty}$, there is a $p_{0} >1$, such that $\omega\in A_{p_{0}}$.
We can choose $\delta>0$ small enough, so that  $p/\delta>p_{0}$. So $\omega\in A_{p/\delta}$.
Then by the definition of $M_{\delta}$ and the weighted $L^{p/\delta}-$boundedness of the Hardy-Littlewood maximal operator $M$, we have
\begin{eqnarray}  \label{equ:3.29}
\begin{split}
\dint_{\mathbf{R}^{n}}[M_{\delta}(\mu_{\Omega}(f))(x)]^{p}\omega(x)\mathrm{d}x &=\dint_{\mathbf{R}^{n}}[M(|\mu_{\Omega}(f)|^{\delta})(x)]^{p/\delta}\omega(x)\mathrm{d}x \\
&\leq \dint_{\mathbf{R}^{n}}|\mu_{\Omega}(f)(x)|^{p}\omega(x)\mathrm{d}x.
\end{split}
\end{eqnarray}

This, together with (\ref{equ:3.10}), Theorem \ref{thm.1} and the fact
$M(f)\leq CM_{L(\log L)^{1/s}}(f)$ for any $s>0$, gives
\begin{eqnarray*}
&\|\mu_{\Omega, b_{1}}(f)\|_{L^{p}(\omega)}&\leq C\|b_{1}\|_{\Osc_{\exp L ^{r_{1}}}}\Big(\|M_{L(\log L)^{1/r_{1}}}(f)\|_{L^{p}(\omega)} +\|\mu_{\Omega}(f)\|_{L^{p}(\omega)}\Big)\\
&&\leq C\|b_{1}\|_{\Osc_{\exp L ^{r_{1}}}}\Big(\|M_{L(\log L)^{1/r_{1}}}(f)\|_{L^{p}(\omega)}
+\|M(f)\|_{L^{p}(\omega)}\Big)\\
&&\leq C\|b_{1}\|_{\Osc_{\exp L ^{r_{1}}}}\|M_{L(\log L)^{1/r_{1}}}(f)\|_{L^{p}(\omega)}.
\end{eqnarray*}

Now, suppose that the theorem is true for $1, 2, \cdots, m-1$ and
let us prove it for $m$.  Recall that, if
$\sigma=\{\sigma(1),\cdots,\sigma(j)\}(1\leq j\leq m)$
 and the corresponding satisfies $1/r_{\sigma}=1/r_{\sigma(1)}+\cdots+1/r_{\sigma(j)}$, then $\sigma'=\{1,\cdots,m\}\setminus\sigma$ and the
corresponding $r_{\sigma'}$  satisfying
$1/r_{\sigma'}=1/r-1/r_{\sigma}$.  Reasoning as in
 (\ref{equ:3.29}), for
$\theta>0$ small enough, we have
\begin{equation}
 \label{equ:3.36}
 \int_{\mathbf{R}^{n}}[M_{\theta}(\mu_{\Omega, \vec{b}_{\sigma'}}(f))(x)]^{p}\omega(x)\mathrm{d}x
\leq C\int_{\mathbf{R}^{n}}|\mu_{\Omega, \vec{b}_{\sigma'}}(f)(x)|^{p}\omega(x)\mathrm{d}x,
\end{equation}

The same argument as used above and the induction hypothesis give us that
\begin{eqnarray*}
&&\|\mu_{\Omega, \vec{b}}(f)\|_{L^{p}(\omega)}\leq \|M_{\delta} (\mu_{\Omega, \vec{b}}(f))\|_{L^{p}(\omega)}
\leq C\|M_{\delta}^{\sharp}(\mu_{\Omega, \vec{b}}(f))\|_{L^{p}(\omega)}\\
&&\leq C\|\vec{b}\|\|M_{L(\log L)^{1/r}}(f)\|_{L^{p}(\omega)} +C\sum_{j=1}^{m}\sum_{\sigma\in\mathscr{C}_{j}^{m}}
\|\vec{b}_{\sigma}\|_{\Osc_{\exp L ^{r_{\sigma}}}}\|M_{\varepsilon}(\mu_{\Omega, \vec{b}_{\sigma'}}(f))\|_{L^{p}(\omega)} \\
&&\leq C\|\vec{b}\|\|M_{L(\log L)^{1/r}}(f)\|_{L^{p}(\omega)} +C\sum_{j=1}^{m}\sum_{\sigma\in\mathscr{C}_{j}^{m}}
\|\vec{b}_{\sigma}\|_{\Osc_{\exp L ^{r_{\sigma}}}}\|\mu_{\Omega, \vec{b}_{\sigma'}}(f)\|_{L^{p}(\omega)} \\
&&\leq C\|\vec{b}\|\|M_{L(\log L)^{1/r}}(f)\|_{L^{p}(\omega)} \\
&&\hskip2em  +C\sum_{j=1}^{m}\sum_{\sigma\in\mathscr{C}_{j}^{m}} \|\vec{b}_{\sigma}\|_{\Osc_{\exp L ^{r_{\sigma}}}}\|\vec{b}_{\sigma'}\|_{\Osc_{\exp L ^{r_{\sigma'}}}}\|M_{L(\log L)^{1/r_{\sigma'}}}(f)\|_{L^{p}(\omega)} \\
&&\leq C\|\vec{b}\|\|M_{L(\log L)^{1/r}}(f)\|_{L^{p}(\omega)},
\end{eqnarray*}
where the fourth inequality follows from  (\ref{equ:3.36}) and the
last one follows from the fact that $M_{L(\log L)^{1/r_{\sigma'}}}(f)\leq M_{L(\log L)^{1/r}}(f)$.

To finish the proof of this special case of Theorem \ref{thm.2}, we need
to check (\ref{equ:3.35}). From (\ref{equ:3.36}), it suffices to
prove
\begin{equation} \label{equ:3.37}
\int_{\mathbf{R}^{n}}|\mu_{\Omega, \vec{b}}(f)(x)|^{p}\omega(x)\mathrm{d}x<\infty, \quad \forall~ 0<p<\infty.
\end{equation}
whenever the weight $\omega$ and the functions $b_{j}~(1\leq j\leq
m)$  are all bounded.

Assume that $\mbox{supp}f\subset B=B(0,R)$ for some $R>0$ and write
\begin{eqnarray*}
\int_{\mathbf{R}^{n}}|\mu_{\Omega,\vec{b}}(f)(x)|^{p}\omega(x)\mathrm{d}x =\int_{2B}|\mu_{\Omega,
\vec{b}}(f)(x)|^{p}\omega(x)\mathrm{d}x+
\int_{(2B)^{c}}|\mu_{\Omega, \vec{b}}(f)(x)|^{p}\omega(x)\mathrm{d}x
=I+II,
\end{eqnarray*}

Noting that $\omega$ and $b_{j}$ are all bounded, by the H\"older
inequality, the induction hypothesis and the fact $M_{L(\log
L)^{m}}\sim M^{m+1}$, $L^{p/\delta}-$boundedness of $M$ , there is
\begin{eqnarray*}
\int_{2B}|b_{\sigma}(x)|^{p} |\mu_{\Omega,b_{\sigma'}}(f)(x)|^{p}\omega(x)\mathrm{d}x \leq C_{\omega}\|b_{\sigma}\|_{L^{\infty}(\mathbf{R}^{n})}^{p}|B|^{1-\delta}\|b_{\sigma'}f\|_{L^{p/\delta}(\mathbf{R}^{n})}^{p}<\infty.
\end{eqnarray*}
This and the definition of $\mu_{\Omega, \vec{b}}(f)$ give us that
\begin{equation} \label{equ:3.38}
I\leq
C\sum_{j=1}^{m}\sum_{\sigma\in\mathscr{C}_{j}^{m}}\int_{2B}|b_{\sigma}(x)|^{p}
|\mu_{\Omega, b_{\sigma'}}(f)(x)|^{p}\omega(x)\mathrm{d}x<\infty.
\end{equation}

To deal with II , we first estimate $\mu_{\Omega, \vec{b}}(f)(x)$ for $x\in(2B)^{c}$. $|x|/2\leq|x-y|\leq 3|x|/2$ when $x\in(2B)^{c}$ and $y\in B$. Noting that $\Omega\in L^{\infty}(S^{n-1})$,
 $\omega$ and $b_{j}$ are bounded functions and $|x|\sim|x-y|$ when $x\in(2B)^{c}$ and $y\in B$,
 there is a constant
$C_{\Omega,\vec{b},\omega}$, depending on the $L^{\infty}$ -norm of
$\Omega$, $b_{j}$ and $\omega$, such that
\begin{equation} \label{equ:3.39}\begin{split}
\mu_{\Omega, \vec{b}}(f)(x)&\leq
C\|\Omega\|_{L^{\infty}(S^{n-1})}\|\vec{b}\|_{L^{\infty}(\mathbf{R}^{n})}\bigg(\dint_{0}^{\infty}\Big|\int_{|x-y|\leq
t}\dfrac{|f(y)|}{|x-y|^{n-1}}\mathrm{d}y\Big|^{2}\dfrac{\mathrm{d}t}{t^{3}}\bigg)^{1/2}\\
&\leq
C_{\Omega,\vec{b},\omega}\dint_{\mathbf{R}^{n}}\dfrac{|f(y)|}{|x-y|^{n-1}}
\bigg(\int_{|x-y|\leq t}\dfrac{\mathrm{d}t}{t^{3}}\bigg)^{1/2}\mathrm{d}y\\
&\leq
C_{\Omega,\vec{b},\omega}\dint_{\mathbf{R}^{n}}\dfrac{|f(y)|}{|x-y|^{n}}\mathrm{d}y
\leq
C_{\Omega,\vec{b},\omega}\dfrac{1}{|2B|}\int_{\mathbf{R}^{n}}|f(y)|\mathrm{d}y\leq
C_{\Omega,\vec{b},\omega}M(f)(x).\end{split}
\end{equation}

By (\ref{equ:3.39}) and the fact that $M(f)(x)\leq CM_{L(\log
L)^{1/r}}(f)(x)$, it follows from (\ref{equ:3.34}) that
$$II\leq C_{\Omega,\vec{b},\omega}\int_{(2B)^{c}}[M_{L(\log L)^{1/r}}(f)(x)]^{p}\omega(x)\mathrm{d}x<\infty.$$

This together with (\ref{equ:3.38}) shows that (\ref{equ:3.37}) is
true when $\omega$ and $b_{j}$  are bounded functions, so does
(\ref{equ:3.35}). And then Theorem \ref{thm.2} is proven for this special case.

Case II\quad    For unbounded $\omega$ and $b_{j}$, we will truncate the
weight $\omega$ and the functions $b_{j}(j=1,\cdots,m)$ as follows.
Let $N$  be a positive integer, denote by
$\omega_{N}=\inf\{\omega,N\}$ and by
$\vec{b}^{N}=(b_{1}^{N},\cdots,b_{m}^{N})$, where $b_{j}^{N}$  is
defined by
$$
 b_{j}^{N}(x)=\left\{
\begin{array}{lll} 
N, &\hbox{when}~ b_{j}(x)>N,   \\
  b_{j}(x), &\hbox{when}~  |b_{j}(x)|\leq N, \\
  -N, &\hbox{when}~  b_{j}(x)<-N.
\end{array} \right.
$$
By Lemma 2.4 in \cite{PT}, there is a positive constant $C$
independent of $N$ such that
\begin{equation} \label{equ:3.40}
\|b_{j}^{N}\|_{\Osc_{\exp L ^{r_{j}}}}\leq \|b_{j}\|_{\Osc_{\exp L ^{r_{j}}}}.
\end{equation}

 Applying (\ref{equ:3.7}) for $\vec{b}^{N}$ and $\omega_{N}$, and using (\ref{equ:3.40}),  we have
\begin{equation} \label{equ:3.41}
\int_{\mathbf{R}^{n}}|\mu_{\Omega,\vec{b}^{N}}(f)(x)|^{p}\omega_{N}(x)\mathrm{d}x \leq
C\|\vec{b}\|^{p}\int_{\mathbf{R}^{n}}[M_{L(\log
L)^{1/r}}(f)(x)]^{p}\omega(x)\mathrm{d}x.
\end{equation}

Next, taking into account the fact that $f$  has compact support, we
deduce that $b_{j}^{N}$  converges to $b_{j}$   and
$b_{\sigma(1)}^{N}\cdots b_{\sigma(j)}^{N}f$ converges to
$b_{\sigma(1)}\cdots b_{\sigma(j)}f$ in any space $L^{p}$ for $p>1$
as $N\rightarrow\infty$.  Recalling the $L^{p}$-boundedness of
$\mu_{\Omega}$, we claim that, at least for a subsequence,
$\{|\mu_{\Omega, \vec{b}^{N}}(f)(x)|^{p}\omega_{N}(x)\}_{N=1}^{\infty}$ converges pointwise almost everywhere
to $|\mu_{\Omega, \vec{b}}(f)(x)|^{p}\omega(x)$ as
$N\rightarrow\infty$.

This fact, together with (\ref{equ:3.41}) and Fatou's lemma, finishes the proof of Theorem \ref{thm.2}.
\end{proof}

\section{Proof of Theorem \ref{thm.3} }

The idea of the proof of Theorem \ref{thm.3} follows that of Theorem 1.5
in \cite{PT}.  We first prove the following lemma.

\begin{lemma} \label{lem.5} \ \ Let $\omega\in A_{\infty}$, $\Phi(t)=t\log^{1/r}(\mathrm{e}+t)$, $\vec{b}, r$, and $r_{j}$ be the same as in Theorem \ref{thm.3}. Then for $\rho>2$, $\Omega\in
L^{\infty}(S^{n-1})$ is homogeneous of degree zero and satisfies (\ref{equ:3.1}) and (\ref{equ:3.4}),there exists a positive constant $C$ such that
\begin{eqnarray} \label{equ:3.42}
\begin{split}
 \sup\limits_{t>0}\dfrac{\omega(\{y\in\mathbf{R}^{n}:
M^{\sharp}_{\delta}(\mu_{\Omega, \vec{b}}(f))(y)>t\})}{\Phi(1/t)} \leq
C\sup\limits_{t>0}\dfrac{\omega(\{y\in\mathbf{R}^{n}:
M_{\Phi}(\|\vec{b}\|f)(y)>t\})}{\Phi(t)}%
\end{split}
\end{eqnarray}
for all bounded functions $f$  with compact support and all $0<\delta<1$.
\end{lemma}

\begin{proof}[\bf\indent Proof]\quad To use Lemma 2.1 (a), we first check that
\begin{equation} \label{equ:3.43}
\sup_{t>0}\frac{1}{\Phi(1/t)}\omega(\{x\in\mathbf{R}^{n}:
M_{\varepsilon}(\mu_{\Omega, \vec{b}}(f))(x)>t\})<\infty
\end{equation}
for all bounded functions $f$ with compact support and all $\delta$
with $0<\delta<1$.

We only prove (\ref{equ:3.43}) for the special case where $\omega$
 and $b_{j}$ are bounded functions. For the general case,  we consider
the truncations of $\omega$  and $\vec{b}$ as in the proof of
Theorem 1.2, by a limit discussion, this time, we take into account
the weak (1,1) boundedness of $\mu_{\Omega}$ that gives the
convergence in measure.  Then we can obtain (\ref{equ:3.43}) for all
$\omega$ and $\vec{b}$ with the hypotheses of Lemma 4.1, we omit the
details.

Assume that $\mbox{supp}f\subset B=B(0,R)$, Then, for any $0<\varepsilon<1$
\begin{eqnarray} \label{equ:3.44}
\begin{split}
&\sup\limits_{t>0}\dfrac{\omega(\{x\in\mathbf{R}^{n}: M_{\varepsilon}(\mu_{\Omega, \vec{b}}(f))(x)>t\})}{\Phi(1/t)} \leq C_{\varepsilon}\sup\limits_{t>0}\dfrac{\omega(\{x\in\mathbf{R}^{n}: M_{\varepsilon}(\chi_{2B}\mu_{\Omega, \vec{b}}(f))(x)>t/2\})}{\Phi(1/t)} \\
&\quad  + C_{\varepsilon} \sup\limits_{t>0}\dfrac{1}{\Phi(1/t)}\omega(\{x\in\mathbf{R}^{n}:
M_{\varepsilon}(\chi_{(2B)^{c}}\mu_{\Omega, \vec{b}}(f))(x)>t/2\})=C_{\varepsilon}(I+II),
\end{split}
\end{eqnarray}
where $C_{\varepsilon}$  is a positive constant depending on $\varepsilon$.

For I, making use of the weak (1,1) boundedness of $M$  and $[\Phi(1/t)]^{-1}\leq
Ct$, and noting that $\omega$ and $b_{j}$ are all bounded. Then there is a positive constant $C_{\omega}$, depending on $\omega$, such that
\begin{eqnarray*}
& I &\leq C_{\omega}\sup_{t>0}t|\{x\in\mathbf{R}^{n}: M_{\varepsilon}(\chi_{2B}\mu_{\Omega, \vec{b}}(f))(x)>t/2\}|\\
&&\leq C_{\omega}\int_{2B} |\mu_{\Omega, \vec{b}}(f)(x)| \mathrm{d}x\leq C_{\omega}|B|^{1/2}\bigg(\int_{2B}|\mu_{\Omega, \vec{b}}(f)(x)|^{2}\mathrm{d}x\bigg)^{1/2}<\infty,
\end{eqnarray*}
where the last step follows as (\ref{equ:3.38}).

Recall the fact that $(M(f))^{\varepsilon}\in
A_{1}$ for $0<\varepsilon<1$ and $f$ locally integrable, then
$$M_{\varepsilon}(M(f))(x)=[M(|M(f)|^{\varepsilon})(x)]^{1/\varepsilon}\leq CM(f)(x).$$

Noting that $\omega$ is bounded, it follows from (\ref{equ:3.39}) and the weak (1,1) boundedness of $M$ that
\begin{eqnarray*}&II&\leq C_{\omega}\sup_{t>0}t \cdot \omega(\{x\in\mathbf{R}^{n}:
M_{\varepsilon}(M(f))(x)>Ct\})\\
&& \leq C_{\omega}\sup_{t>0}t \cdot \omega(\{x\in\mathbf{R}^{n}: M(f)(x)>Ct\})\\
&& \leq C_{\omega} \int_{\mathbf{R}^{n}}|f(x)|\mathrm{d}x<\infty.
\end{eqnarray*}

Combining (\ref{equ:3.44}) and the estimates for I and II , we have (\ref{equ:3.43}).

Now, let us turn to proving (\ref{equ:3.42}) by induction.  For $\vec{b}\in \Osc_{\exp L^{r}}$, write
$\tilde{b}=\vec{b}/\|\vec{b}\|$, then $\|\tilde{b}\|=1$, and $\mu_{\Omega, \vec{b}}(f)/\|\vec{b}\|=\mu_{\Omega, \vec{b}/\|\vec{b}\|}(f)=\mu_{\Omega, \tilde{b}}(f)$. So we can
assume that $\|\vec{b}\|=1$. For $m=1$, we understand $\vec{b}=b$, $\|\vec{b}\|=\|b\|_{\Osc_{\exp L^{r}}}=1$, $\mu_{\Omega, \vec{b}}(f)=\mu_{\Omega, b}(f)$. Therefore, to prove (\ref{equ:3.42}), it suffices to prove
\begin{eqnarray} \label{equ:3.45}
\begin{split}
\sup\limits_{t>0}\dfrac{\omega(\{y\in\mathbf{R}^{n}: M^{\sharp}_{\delta}(\mu_{\Omega, b}(f))(y)>t\})}{\Phi(1/t)} \leq C\sup\limits_{t>0}\dfrac{\omega(\{y\in\mathbf{R}^{n}: M_{L(\log L)^{1/r}}(f)(y)>t\})}{\Phi(1/t)}
\end{split}
\end{eqnarray}
for all bounded functions $f$ with compact support.

Applying Lemma \ref{lem.4} for $m=1$ and any $\varepsilon$ with
$0<\delta<\varepsilon<1$, it is easy to see that the left-hand side
of (\ref{equ:3.45}) is dominated by
\begin{eqnarray*}
\sup\limits_{t>0}\dfrac{\omega(\{y\in\mathbf{R}^{n}:M^{\sharp}_{\delta}(\mu_{\Omega, b}(f))(y)>t\})}{\Phi(1/t)} &\leq&
C\sup_{t>0}\frac{\omega(\{y\in\mathbf{R}^{n}: M_{L(\log L)^{1/r}}(f)(y)>t/2\})}{\Phi(1/t)}\\
&&+C\sup_{t>0}\frac{\omega(\{y\in\mathbf{R}^{n}:M_{\varepsilon}(\mu_{\Omega}(f))(y)>t/2\})}{\Phi(1/t)}.
\end{eqnarray*}

Recall that (\ref{equ:3.43}) is valid and since $[\Phi(1/t)]^{-1}$
is doubling, then by Lemma \ref{lem.1} (a), Lemma \ref{lem.3}  and noting that
$M(f)\leq M_{L(\log L)^{1/r}}(f)$, we have
\begin{eqnarray*}
&& \sup_{t>0}\frac{\omega(\{y\in\mathbf{R}^{n}:M^{\sharp}_{\delta}(\mu_{\Omega, b}(f))(y)>t\})}{\Phi(1/t)} \leq C\sup_{t>0}\frac{\omega(\{y\in\mathbf{R}^{n}: M_{L(\log L)^{1/r}}(f)(y)>t\})}{\Phi(1/t)}\\
&&\qquad + C\sup_{t>0}\frac{\omega(\{y\in\mathbf{R}^{n}: M_{\varepsilon}^{\sharp}(\mu_{\Omega}(f))(y)>t\})}{\Phi(1/t)}\\
&& \leq C\sup_{t>0}\frac{\omega(\{y\in\mathbf{R}^{n}: M_{L(\log L)^{1/r}}(f)(y)>t\})}{\Phi(1/t)} + C\sup_{t>0}\frac{\omega(\{y\in\mathbf{R}^{n}: M(f)(y)>t\})}{\Phi(1/t)} \\
&& \leq C\sup_{t>0}\frac{1}{\Phi(1/t)}\omega(\{y\in\mathbf{R}^{n}: M_{L(\log L)^{1/r}}(f)(y)>t\}).
\end{eqnarray*}
This is (\ref{equ:3.45}), thus, we have proved (\ref{equ:3.42}) for $m=1$.

Now, let us check (\ref{equ:3.42}) for the general case $m\geq 2$.
Suppose that (\ref{equ:3.42}) holds for $m-1$, let us prove it for
$m$. Noting that (\ref{equ:3.43}) is true and recalling the fact
that $[\Phi(1/t)]^{-1}$ is doubling, then by Lemma \ref{lem.3} and \ref{lem.4} for
$\varepsilon$ with $0<\delta<\varepsilon$, Lemma \ref{lem.1} (a) and the
induction hypothesis on (\ref{equ:3.42}), we have
\begin{eqnarray*}
&&\sup_{t>0}\frac{\omega(\{y\in\mathbf{R}^{n}: M^{\sharp}_{\delta}(\mu_{\Omega, \vec{b}}(f))(y)>t\})}{\Phi(1/t)}  \leq C\sup_{t>0}\frac{\omega(\{y\in\mathbf{R}^{n}: M_{\Phi}(f)(y)>t/C_{m}\})}{\Phi(1/t)}\\
&&\qquad +C\sum_{j=1}^{m} \sum_{\sigma\in \mathscr{C}_{j}^{m}}\sup_{t>0}\frac{1}{\Phi(1/t)}\omega(\{y\in\mathbf{R}^{n}: M_{\varepsilon}(\mu_{\Omega, \vec{b}_{\sigma'}}(\|\vec{b}_{\sigma}\|f))(y)>t/C_{m}\})\\
&& \leq C_{m} \sup_{t>0}\frac{1}{\Phi(1/t)}\omega(\{y\in\mathbf{R}^{n}: M_{\Phi}(f)(y)>t\})\\
&&\qquad+C_{m}\sum_{j=1}^{m} \sum_{\sigma\in \mathscr{C}_{j}^{m}}\sup_{t>0}\frac{1}{\Phi(1/t)}\omega(\{y\in\mathbf{R}^{n}:
M_{\varepsilon}^{\sharp}(\mu_{\Omega, \vec{b}_{\sigma'}}(\|\vec{b}_{\sigma}\|f))(y)>t\})\\
&& \leq C_{m}\sup_{t>0}\frac{1}{\Phi(1/t)}\omega(\{y\in\mathbf{R}^{n}: M_{\Phi}(f)(y)>t\})\\
&&\qquad+C_{m}\sum_{j=1}^{m} \sum_{\sigma\in \mathscr{C}_{j}^{m}}\sup_{t>0}\frac{1}{\Phi(1/t)}\omega(\{y\in\mathbf{R}^{n}:
M_{\Phi}(\|\vec{b}_{\sigma'}\|\|\vec{b}_{\sigma}\|f)(y)>t\})\\
&& \leq C_{m}\sup_{t>0}\frac{1}{\Phi(1/t)}\omega(\{y\in\mathbf{R}^{n}: M_{\Phi}(f)(y)>t\})\\
&&\qquad+C_{m}\sum_{j=1}^{m} \sum_{\sigma\in \mathscr{C}_{j}^{m}}\sup_{t>0}\frac{1}{\Phi(1/t)}\omega(\{y\in\mathbf{R}^{n}: M_{\Phi}(f)(y)>t\}),
\end{eqnarray*}
where $\|\vec{b}_{\sigma}\|$ and $\|\vec{b}_{\sigma'}\|$ are as in (\ref{equ:3.15}), and in the last step, we make
use of the fact that $\|\vec{b}_{\sigma'}\|\|\vec{b}_{\sigma}\|=\|\vec{b}\|=1$.

    This concludes (\ref{equ:3.42}) for all $m$, so the proof of Lemma \ref{lem.5} is completed.
\end{proof}

\begin{lemma}\label{lem.6} \ \ Let  $ \omega\in A_{\infty}$, $\Phi(t)=t\log^{1/r}(\mathrm{e}+t)$, $\vec{b}, r$, and $r_{j}$  be the same as in Theorem \ref{thm.3}.  For $\rho>2$, $\Omega\in
L^{\infty}(S^{n-1})$ is homogeneous of degree zero and satisfies (\ref{equ:3.1}) and (\ref{equ:3.4}), there exists a positive constant $C$ such that
$$\sup_{t>0}\frac{\omega(\{y\in\mathbf{R}^{n}: \mu_{\Omega, \vec{b}}(f)(y)>t\})}{\Phi(1/t)} \leq
C\sup_{t>0}\frac{\omega(\{y\in\mathbf{R}^{n}: M_{\Phi}(\|\vec{b}\|f)(y)>t\})}{\Phi(1/t)}$$
for all bounded functions $f$ with compact support.
\end{lemma}

The proof is similar as the proof of Lemma 4.2 in \cite{Z2}, we omit the details here.

To prove Theorem  \ref{thm.3}, we need the following weighted weak-type inequality due to P\'{e}rez and Trujillo-Gonz\'{a}lez\upcite{PT}.

\begin{lemma}\upcite{PT}\label{lem.7}\ \ Let $\omega\in A_{1}$, $\Phi(t)=t\log^{1/r}(\mathrm{e}+t)$. Then there is a positive constant $C$, for any $\lambda>0$ and any locally
integrable function $f$ , such that
$$\omega(\{y\in\mathbf{R}^{n}: M_{\Phi}(f)(y)>\lambda\}) \leq
C\int_{\mathbf{R}^{n}}\Phi\bigg(\frac{|f(y)|}{\lambda}\bigg)\omega(y)\mathrm{d}y.$$
\end{lemma}

 \begin{proof} [\indent\bf Proof of Theorem \ref{thm.3}]  \  By homogeneity of $\vec{b}$, we can assume that $\lambda=\|\vec{b}\|=1$. Then we only need to prove that
$$\omega(\{y\in\mathbf{R}^{n}: \mu_{\Omega, \vec{b}}(f)(y)>1\})
\leq C\int_{\mathbf{R}^{n}}\Phi\big(|f(y)|\big)\omega(y)\mathrm{d}y.
$$

By $\Phi(ab)\leq 2\Phi(a)\Phi(b)$, $a,b\geq 0$ and Lemma \ref{lem.6}, Lemma \ref{lem.7}, we have
\begin{eqnarray*}
&& \omega(\{y\in\mathbf{R}^{n}:\mu_{\Omega, \vec{b}}(f)(y)>1\}) \leq C\sup_{\lambda>0}\frac{1}{\Phi(1/\lambda)}\omega(\{y\in\mathbf{R}^{n}:\mu_{\Omega, \vec{b}}(f)(y)>\lambda\})\\
&&\leq C\sup_{\lambda>0}\frac{\omega(\{y\in\mathbf{R}^{n}:M_{\Phi}(f)(y)>\lambda\})}{\Phi(1/\lambda)}  \leq C\sup_{\lambda>0}\frac{1}{\Phi(1/\lambda)}\int_{\mathbf{R}^{n}}\Phi\bigg(\frac{|f(y)|}{\lambda}\bigg)\omega(y)\mathrm{d}y\\
&&\leq
C\sup_{\lambda>0}\frac{1}{\Phi(1/\lambda)}\int_{\mathbf{R}^{n}}\Phi\big(|f(y)|\big)\Phi(1/\lambda)\omega(y)\mathrm{d}y \leq C\int_{\mathbf{R}^{n}}\Phi\big(|f(y)|\big)\omega(y)\mathrm{d}y.
\end{eqnarray*}

Then we finish the proof of Theorem \ref{thm.3}.
\end{proof}

\begin{acknowledgments}
 The authors cordially  thank the referees for their valuable suggestions and useful comments which have lead to the improvement of this paper. This work was   supported  in part  by the Pre-Research Project of Provincial Key Innovation (No. SY201224),   the NSF  (No. A200913) of Heilongjiang Province and NNSF (No. 11041004 and 11161042)  of China.
\end{acknowledgments}


\bigskip
\address{ Jianglong Wu\\
Department of Mathematics \\
Mudanjiang Normal University \\
Mudanjiang 157011 \\
China } {jl-wu@163.com}
     \address{Qingguo Liu\\
Laboratory for Multiphase Processes \\
University of Nova Gorica \\
Nova Gorica 5000 \\
Slovenia } {liuqingguo1980@gmail.com}

\end{document}